\documentclass{article}

\usepackage{amsmath, amsthm, amsfonts,amssymb}
\usepackage{color}

\usepackage{hyperref}

\usepackage{longtable}

\usepackage{csquotes}

\newtheorem{thm}{Theorem}[section]
\newtheorem{cor}[thm]{Corollary}
\newtheorem{lem}[thm]{Lemma}

\usepackage{mdframed}
\theoremstyle{definition}
\newmdtheoremenv{boxProb}{Problem}
\newmdtheoremenv{boxDef}{Definition}
\newmdtheoremenv{boxCor}{Corollary}
\newmdtheoremenv{boxThm}{Theorem}
\newmdtheoremenv{compjob}{Computational Job}
\newmdtheoremenv{reqi}{Requirement}

\usepackage{graphicx}
\usepackage{subfigure}          

\usepackage{placeins}

\usepackage{tabularx}

\usepackage{lmodern,textcomp}

\usepackage{algorithm}
\usepackage{algpseudocode}

\usepackage{enumitem}

\usepackage{xspace} 
\newcommand\largeparbreak{\par\bigskip}

\newcommand*\tageq{\refstepcounter{equation}\tag{\theequation}}

\newcommand{\inv}{^{-1}\xspace}

\newcommand{\blambda}{{\boldsymbol{\lambda}}\xspace}

\renewcommand{\t}{^\textsf{T}\xspace}

\newcommand{\away}[1]{}

\newcommand{\R}{\mathbb{R}\xspace}

\newcommand{\N}{\mathbb{N}\xspace}

\newcommand{\cO}{\mathcal{O}\xspace}

\newcommand{\bA}{\mathbf{A}\xspace}

\newcommand{\bG}{\mathbf{G}\xspace}

\newcommand{\bZ}{\mathbf{Z}\xspace}

\newcommand{\bb}{\mathbf{b}\xspace}

\newcommand{\bei}[1]{{\mathbf{e}}\xspace}

\newcommand{\bM}{\mathbf{M}\xspace}

\newcommand{\bI}{\mathbf{I}\xspace}

\newcommand{\bu}{\mathbf{u}\xspace}
\newcommand{\bv}{\mathbf{v}\xspace}

\newcommand{\tlambda}{{\tilde{\lambda}}\xspace}

\usepackage{todonotes}

\title{Determination of Positive Definiteness through Shift-and-Invert Iteration in Weakly Polynomial Complexity}

\author{Martin Neuenhofen}

\begin{document}

\maketitle

\begin{abstract}
We propose a numerical method, based on the shift-and-invert power iteration, that answers whether a symmetric matrix is positive definite ("yes") or not ("no").

Our method uses randomization. But, it returns the correct answer with high probability. A thorough proof for the probability is presented. If the method answers "yes", the result is true with a high constant probability. If it answers "no", it provides proof that the matrix is not positive definite.

The method has the following benefits: The cost for a constant probability of success scales logarithmically with the condition number. Further, since essentially consisting of vector iterations, our method is easy to implement.
\end{abstract}


\subsubsection*{Brief Summary} The method is on page 10. The iteration complexity of the method is given on page 11.

\section{Introduction}
This paper proposes a method that determines whether a given real symmetric matrix is positive definite or not.

This problem, also referred to as identification of positive definiteness, is fundamental in matrix analysis. Solving this problem is desirable for the following applications:
\begin{enumerate}
	\item In Quasi-Newton methods, inertia correction schemes are used to safeguard global convergence of the iteration.
	\item One important subclass of Quasi-Newton methods are optimization algorithms. Here, inertia-correction schemes are used to make sure that the step-direction is a descent direction for the objective to be minimized.
	\item In Successive Quadratic Programming (SQP) methods, the Hessian is typically modified into a semi-definite form, such that the step-direction can be computed in efficient time-complexity.
\end{enumerate}

\largeparbreak

Two naive approaches are used commonly for the determination of positive definiteness of a matrix.

The first approach works by using Cholesky decomposition. If the matrix is positive definite, then the decomposition is stable and will hence succeed. A certificate of positive definiteness is provided through the Cholesky factor, that can be thought of as a square-root of the matrix. Clearly, it only exists when all the eigenvalues of the matrix are strictly positive.

The second approach works by computing the smallest eigenvalue of the matrix\footnote{Readers beware: In the general case, the smallest eigenvalue of a symmetric matrix is not identical to the eigenvalue of smallest absolute value.}. The matrix is positive definite if and only if this eigenvalue is positive. Yet, when using an iterative scheme to compute the smallest eigenvalue, there is the risk that the iteration traps into a semi-minimal eigenvalue. For example, when using a Rayleigh quotient iteration with a random starting vector and an initial guess of the eigenvalue, then this iteration is quite likely to converge towards an eigenvalue that is not the closest one to the initial guess!

\largeparbreak

In this work we are interested in the setting where the matrix is large and sparse, and where computations can be performed on massively parallel machines. Discussing the above two approaches with respect to this setting, it is certainly unattractive to perform a Cholesky decomposition. This is because Cholesky has limited potential for massively parallel computations. Vector iterations are more promising, as each iteration can be parallelized easily. In contrast to a Cholesky decomposition, one iteration of vector iteration can be realized in a constant time-complexity, irrespective of the dimension of the system matrix. For banded matrices of small bandwidth, even (shifted) inverse iterations can be computed in almost constant parallel time-complexity \cite{myblocktridiagLGS}. As a result, for the determination of positive definiteness, power-type iteration method appear most suitable. Yet, there are a number of issues that must be overcome.

One issue with power-type iterations is the convergence to the correct eigenvalue. For instance, let us consider the Rayleigh-quotient iteration. This iteration has a rapid rate of convergence, in fact a cubic rate \cite[Chapter~8]{MatrixComputations}. Hence, we find a highly accurate approximation for an eigenvalue of the matrix in just a few iterations. But, this eigenvalue is likely to be non-extremal. Thus, apparently the shift-and-invert acceleration technique cannot be used without safeguards when aiming for a reliable identification algorithm. 

Another issue is the speed of convergence. For instance, let us consider the power iteration. Assuming the non-pathologic case, this iteration will always converge to the eigenvalue with largest absolute value. But it may converge very slowly. It could be enhanced with a Lanczos method \cite[Chapter~9]{MatrixComputations}. But in the worst case, both the convergence of Lanczos iteration and power iteration depend strongly on the spectrum.

\largeparbreak

To arrive at a reliable method with only logarithmic iteration complexity on the conditioning, we need a further ingredient.

The key to the idea of our method is the following insight on power iteration applied to symmetric positive definite matrices: While power iteration converges slowly to the largest eigenvalue, it converges rapidly to an approximate eigenvalue that has a relative error of $40\%$ \ --- no matter what spectrum the matrix has. As we are going to show, this is entirely sufficient to construct an efficient method.

\largeparbreak

This paper is organized as follows. Section 2 gives preliminaries. In section 3 we present the method. In section 4 we give conclusions.

\section{Preliminaries}\label{sec:Prelim}
\paragraph{Problem statement}
We consider the situation where we are given a regular real symmetric matrix
\begin{align*}
\bA = \bA\t \in \R^{N \times N}\,,\quad\quad N \in \N\,.
\end{align*}
In the domain of our interest, $N$ is very large. We shall determine with a reliability of $99.9999\%$ whether $\bA$ is positive definite or not. The method only accesses the matrix only via products with $\bA$ and products with shifted inverses of $\bA$.

\paragraph{Notation}
We introduce some notation. We will discuss properties of the spectrum of various symmetric matrices. We write $\sigma(\bG) \subset \R$ for the spectrum of a symmetric matrix $\bG \in \R^{N \times N}$. We indicate the eigenvalues in descending order
\begin{align*}
	\lambda_1(\bA) \geq \lambda_2(\bA) \geq \dots \geq \lambda_N(\bA)\,.
\end{align*}
When clear from the context, we may drop the argument $\bA$.

We consider four particular eigenvalues of the spectrum:
\begin{itemize}
	\item The largest eigenvalue: $\lambda_{\text{max}}(\bG) := \lambda_1$
	\item The smallest eigenvalue: $\lambda_{\text{min}}(\bG) := \lambda_N$
	\item The eigenvalue with largest absolute value: $|\lambda_{\text{absmax}}(\bG)| := \max( |\lambda_1,\lambda_N| )$
	\item The eigenvalue with smallest absolute value: $|\lambda_{\text{absmin}}(\bG)|$.
\end{itemize}
The latter value has no simple expression in general. For positive definite matrices it holds 
$$	|\lambda_{\text{absmax}}|=\lambda_N\,. 	$$

We use the writing $\langle \bu , \bv \rangle$ to denote the scalar-product between two vectors $\bu,\bv$ of equal dimension. We use $\|\cdot\|_p$ for the $p$-norm of square matrices and vectors.
We write the condition number of a regular matrix $\bG \in \R^{N \times N}$ as
$$ 		\kappa(\bG) := \|\bG\|_2 \cdot \|\bG\inv\|_2\,. 	$$

We use the notation of a random oracle $B_N$, that returns uniformly distributed random samples from the unit-circle of $\R^N$. We write $\bb:=B_N$ to instantiate a vector $\bb\in\R^N$ based on this distribution.
The following result from \cite{RandomSphereVector} describes how $B_N$ can be implemented in a practical way.
\begin{lem}
	Let $\xi_1,\xi_2,...,\xi_n\in\R$ be random samples of a normal distribution. Then a vector $\bb$ as described above can be generated as follows:
	\begin{align*}
	\bb := \frac{1}{\sqrt{\sum_{j=1}^n \xi^2_j}} \cdot \begin{pmatrix}
	\xi_1\\
	\xi_2\\
	\vdots\\
	\xi_N
	\end{pmatrix}
	\end{align*}
\end{lem}

\subsection{Power iteration}

A vital scheme for our proposed method is the power iteration. In this paper, we use the power iteration to compute crude numerical estimates to the largest eigenvalue of various symmetric positive definite matrices. The scheme is defined in Algorithm~\ref{algo:PowerIter}.

\begin{algorithm}
	\caption{Power Iteration}
	\label{algo:PowerIter}
	\begin{algorithmic}[1]
		\Procedure{PowerIteration}{$\bG,\bb,k$}
			\State $\bv := \bb$
			\For{$j=0,1,2,\dots,k$}
				\State $\bb := 1/\|\bv\| \cdot \bv$
				\State $\bv := \bG \cdot \bb$
				\State $\tilde{\lambda} := \langle \bv,\bb\rangle / \langle \bb,\bb \rangle$
			\EndFor
			\State \Return $\tilde{\lambda},\,\bb$
		\EndProcedure
	\end{algorithmic}
\end{algorithm}

Two results for this use case are presented. The first is taken from \cite[Theorem~3.1(a)]{ErrorSymEigvaluePower}.
\begin{thm}[Power iteration]\label{thm:PowerIter}
	Given a symmetric positive definite definite matrix $\bG \in \R^{N \times N}$. Choose $\varepsilon,\delta \in (0,1)$ and generate $\bb:= B_N$. Then, after
	\begin{align}
		k := \left\lceil\,\frac{\log(N/\delta^2)}{2 \cdot \varepsilon}\,\right\rceil\label{eqn:FormulaK}
	\end{align}
	iterations of Algorithm~\ref{algo:PowerIter}, a value $\tlambda$ is returned that satisfies
	\begin{align*}
		\left|\,\frac{\lambda_{\text{max}}(\bG) - \tlambda}{\lambda_{\text{max}}(\bG)}\,\right| \leq \varepsilon
	\end{align*}
	with a probability of at least $1-\delta$.
\end{thm}

In what follows we prepare a technical result that will be useful for our method. It can be thought of as a bisection method for determining the largest eigenvalue of a matrix.

\subsection{Bisection of bound of spectrum}\label{sec:PrelimBisec}
The question whether a symmetric matrix $\bA$ is positive definite or not can be reduced into the problem of identifying, whether the largest eigenvalue of an auxiliary matrix $\bM$ with $\sigma(\bM) \subset(0,2)$ is smaller than $1$. We will approach this problem with a bisection algorithm, that iteratively shrinks an interval in which the largest eigenvalue of $\bM$ lives. Our bisection algorithm is one-sided, i.e. it only shrinks the interval from the right.

\largeparbreak

The bisection makes the initial assumption that $\lambda_{\text{max}}(\bM)<\mu$ holds for some given value $\mu \in (1,2]$.
We define
$$ 	\overline{\mu} := \frac{1+\mu}{2}\,,	$$
which is at half the distance of $\mu$ to $1$.

We wish to determine whether $\lambda_{\text{max}}(\bM)<\overline{\mu}$ holds with probability $1-\delta$, for some very small chosen $\delta\in(0,1)$\,. To this end, let us consider the following matrix
\begin{align}
	\bZ := (\mu-1)\cdot(\mu\cdot\bI - \bM)\inv\,.
\end{align}
We notice that $\bZ$ is well-defined and positive definite because $\mu$ exceeds the largest eigenvalue of $\bM$.

We chose $\varepsilon:=0.4$ and $\delta$ small as described above. We compute $k$ from $\varepsilon,\delta$, as defined in \eqref{eqn:FormulaK}. After $k$ iterations of power iteration we obtain a returned value $\tlambda$.

All we want to make sure at this point is that either of the following two holds:
\begin{itemize}
	\item Case 1: Our original purpose is to determine whether $\bM$ has any eigenvalue that is larger or equal to $1$. Thus, if $\tlambda\geq 1$, then since $\tlambda \in 1-\sigma(\bM)$, we can be sure that $\lambda_{\text{max}}(\bM)\geq 1$, and hence $\bA$ is indefinite. In this case, we have solved the problem and can terminate.
	\item Case 2: Our purpose for this iteration of the bisection method is to make sure that $\lambda_{\text{max}}(\bM)<\overline{\mu}$. We will show below that if $\tlambda<1$, then this bound holds true with probability $1-\delta$.
\end{itemize}
In summary, either case 1 or case 2 occurs. In case 1, which can be thought of as the lucky case, we terminate the bisection because we know that $\lambda_{\text{max}}(\bM)\geq 1$. In case 2, it holds $\sigma(\bM)\subset(0,\overline{\mu})$ with high probability, which is the initial assumption for bisection with $\overline{\mu}$. Hence, in case 2 we iteratively perform a subsequent bisection step with the above procedure, where we use the bisection update $\mu := \overline{\mu}$. Details of the exact iteration, and when it terminates, will be given in Section~\ref{sec:TheMethod}.

Eventually, we now prove the above claims on cases 1 and 2.
\begin{cor}[Shift-and-invert bisection iteration]
	Let $\mu\in(1,2]$ and $\bM\in\R^{N \times N}$ be symmetric positive definite with $\sigma(\bM)\subset(0,\mu)$. Define
	\begin{align*}
		\overline{\mu} &:= \frac{1+\mu}{2}\,,\\
		\bZ &:= (\mu-1)\cdot(\mu\cdot\bI-\bM)\inv\,.
	\end{align*}
	Choose $\varepsilon:=0.4$, $\delta \in (0,1)$, and define $k$ as in \eqref{eqn:FormulaK}. Generate $\bb:=B_N$. Perform $k$ power iterations with matrix $\bZ$ and vector $\bb$, and obtain output $\tlambda$.
	
	Then the following hold:
	\begin{enumerate}
		\item If $\tlambda\geq 1$ then $\lambda_{\text{max}}(\bM)\geq 1$.
		\item If $\tlambda< 1$ then $\lambda_{\text{max}}(\bM)<\overline{\mu}$ holds with probability $1-\delta$.
	\end{enumerate}
\end{cor}
\noindent
\underline{Proof:}
Let us start with some notes on $\bZ$. The eigenvalues of $\bZ$ and $\bM$ are related as
\begin{align*}
	\lambda_j(\bZ) = \frac{\mu-1}{\mu-\lambda_j(\bM)}\,.
\end{align*}
Notice from the equation, that if $\bZ$ has an eigenvalue $\lambda_j(\bZ)\geq 1$, then this means $\bM$ has an eigenvalue $\lambda_j(\bM)\geq 1$, too. Notice further that if $\lambda_{\text{max}}(\bZ)< 2$, then this implies $\lambda_{\text{max}}(\bM)< \overline{\mu}$. This helps us to show the propositions:

The first proposition can be shown from the following fact: The Rayleigh-quotient, as which $\tlambda$ is computed, is bounded by the spectrum of $\bZ$ \cite{ErrorSymEigvaluePower}. Thus,
$$ 	\lambda_{\text{max}}(\bZ) \geq \tlambda\,. 	$$
Hence, if $\tlambda\geq 1$ then $\lambda_{\text{max}}(\bZ)\geq 1$, and therefore $\lambda_{\text{max}}(\bM)\geq 1$.

The second proposition follows from the probability in which $\tlambda$ is a $40\%$ relatively accurate approximation of an eigenvalue of $\bZ$. If $\tlambda<1$, then the bound
\begin{align*}
	\left|\,\frac{\lambda_{\text{max}}(\bZ)-\tlambda}{\lambda_{\text{max}}(\bZ)}\,\right| \leq 0.4
\end{align*}
holds with probability $1-\delta$. Hence, $\tlambda \geq 0.6 \cdot \lambda_{\text{max}}(\bZ)$ with probability $1-\delta$. In consequence of this, if $\tlambda < 1$ then
$$ 	\lambda_{\text{max}}(\bZ) < \frac{1}{0.6} \cdot 1 < 2\,. 	$$
As we discussed above, this implies $\lambda_{\text{max}}(\bM)<\overline{\mu}$.\qed

\section{The Proposed Method}\label{sec:TheMethod}

\paragraph{Overview}
Our method works in three steps.
\begin{enumerate}
	\item We start by estimating an interval $(\check{\mu},\hat{\mu})\subset(0,\infty)$, that contains with high probability the spectrum of $\bA$. Our estimate we compute with power iterations on the matrices $\bA$ and $\bA\inv$. We conduct these iterations in a way such that we satisfy the bound $\hat{\mu}/\check{\mu} \in \cO(\,\kappa(\bA)\,)$ with high probability.
	\item We then construct a matrix $\bM:= \bI - 1/\hat{\mu} \cdot \bA$\,. This matrix is symmetric positive definite with high probability. If $\bM$ is positive definite, then the matrix $\bA$ is positive definite if and only if $\lambda_{\text{max}}(\bM)<1$.
	\item We perform iteratively a bisection algorithm on the matrix $\bM$ to determine whether the largest eigenvalue is $< 1$. We will show that $\cO(\log(\kappa(\bA)))$ bisection iterations need to be performed for this purpose.
\end{enumerate}

This section is organized as follows. We dedicate one subsection on each step of the method. Eventually we state the whole algorithm with particular values for the parameters. For these parameters we then prove the probability of correctness of our method, and the time-complexity.

\subsection{Accurate Estimation of Bounds of the Spectrum}
We choose $\delta \in (0,1)$ suitably small, $\varepsilon:= 0.4$, and compute $k$ from \eqref{eqn:FormulaK}. We compute $\check{\mu}$ and $\hat{\mu}$ as follows:
\begin{algorithmic}[]
	\State $\tlambda_N^{-1} := $\textsc{PowerIteration}($\bA^{-2},B_N,k$)
	\State $\check{\mu} := \sqrt{0.5 \cdot \tlambda_N}$
	\State $\tlambda_1 := $\textsc{PowerIteration}($\bA^2,B_N,k$)
	\State $\hat{\mu} := \sqrt{1.5 \cdot \tlambda_1}$
\end{algorithmic}
In this we used the trick that the square of a symmetric matrix is positive definite and thus Theorem~\ref{thm:PowerIter} applies.

According to the theorem, $\sigma(\bA) \in (\check{\mu},\hat{\mu})$ holds with high probability. We explain this. Due to Theorem~\ref{thm:PowerIter}, we have computed approximations $\tlambda_1,\tlambda_N$ of relative accuracy $0.4$ with respect to the eigenvalues $\lambda_{\text{absmax}},\lambda_{\text{absmin}}$ of $\bA^2$. Each relative error holds with a probability of $1-\delta$. Hence, both error bounds hold with a probability of $(1-\delta)^2 > 1-2 \cdot \delta$. Further, since we computed $\check{\mu}^2,\hat{\mu}^2$ by decreasing/ increasing $\tlambda_1,\tlambda_N$ by a value that exceeds the upper bound of their relative error, it holds with the above probability that
\begin{align*}
	& 				&	\sigma(\bA^2) &\in (\check{\mu}^2,\hat{\mu}^2)\\
	&\Rightarrow	&	\sigma(\bA)   &\in (\check{\mu},\hat{\mu}) \quad \text{or $\bA$ is not positive definite}
\end{align*} 
for the open interval.

It is important for our complexity analysis, that
\begin{align*}
&	&\frac{\hat{\mu}^2}{\check{\mu}^2} &\leq \frac{1.5}{0.5} \cdot \frac{\tlambda_1}{\tlambda_N} \leq \frac{1.5}{0.5} \cdot \frac{\frac{1}{1-0.4}}{\frac{1}{1+0.4}} \cdot \frac{\lambda_1(\bA^2)}{\lambda_N(\bA^2)} \leq 7 \cdot \kappa(\bA^2) = 7 \cdot \kappa(\bA)^2\\
&\Rightarrow & \frac{\hat{\mu}}{\check{\mu}} &\leq \sqrt{7} \cdot \kappa(\bA)\tageq
\end{align*}
holds. The bound is satisfied with the above high probability because $\hat{\mu}^2,\check{\mu}^2$ are computed as scaled versions of $\tlambda_1,\tlambda_N$, whose relative errors in turn with respect to $\lambda_1,\lambda_N$ are bounded by $0.4$\,.

\subsection{Construction of the auxiliary matrix}
We define an important matrix, that we call \textit{auxiliary matrix}:
\begin{align}
	\bM := \bI - \frac{1}{\hat{\mu}} \cdot \bA
\end{align}
The spectrum of this matrix is bounded as
\begin{align}
\sigma(\bM) \subset \left(\,0\,,\,1-\frac{\check{\mu}}{\hat{\mu}}\,\right]\cup\left[\,1+\frac{\check{\mu}}{\hat{\mu}}\,,\,2\,\right)\,,\label{eqn:SpectrumM}
\end{align}
as can be verified from the bounds on $\sigma(\bA)$ and the construction of $\bM$ from $\bA$.

We make two important observations from \eqref{eqn:SpectrumM}: If $\bA$ is positive definite, then the largest eigenvalue of $\bM$ is bounded above as
\begin{align}
\lambda_{\text{max}}(\bM) \leq 1-\check{\mu}/\hat{\mu}\,. 	\label{eqn:CasePosDef}
\end{align}
Instead, if $\bA$ is not positive definite, then the largest eigenvalue of $\bM$ is bounded below as
\begin{align}
	\lambda_{\text{max}}(\bM) \geq 1+\check{\mu}/\hat{\mu}\,. 	\label{eqn:CaseNotPosDef}
\end{align}

\largeparbreak

\subsection{Shift-and-Invert Probing Iteration}

Define
$$ 	J := \lceil\,\log_2(\hat{\mu}/\check{\mu})\,\rceil	$$
and the following sequence $\lbrace\mu_j\rbrace$\,:
\begin{align*}
	\mu_j := 1+\frac{1}{2^j},\quad\quad j=0,1,2,3,\dots,J\,.
\end{align*}
Further, define the matrices
\begin{align*}
	\bZ_j := (\mu-1) \cdot (\mu_j \cdot \bI-\bM)\inv
\end{align*}
These definitions allow us to state the bisection approach described in Section~\ref{sec:PrelimBisec} in a succinct form:

\largeparbreak

The iterations of the bisection method have the iteration index $j \in \N_0$. The $j$th bisection iteration consists of applying the power iteration for $\bZ_j$. When $\tlambda$, the estimate found by power iteration, is greater or equal to $1$, then we terminate with certificate of indefiniteness of $\bA$. Otherwise, we continue with the next bisection iteration for $j:= j+1$\,.

\paragraph{Termination of the iteration}
Consider the ultimate bisection iteration, i.e. $j=J-1$. At this point we have found out that
$$ 	\lambda_{\text{max}}(\bM)<\mu_J\leq 1+\frac{\check{\mu}}{\hat{\mu}} 	$$
holds with high probability. Combining this result with \eqref{eqn:SpectrumM} yields
\begin{align*}
	\sigma(\bM) \subset \left(\,0\,,\,1-\frac{\check{\mu}}{\hat{\mu}}\,\right]\,.
\end{align*}
with high probability. Hence, we can conclude that $\bA$ is positive definite with high probability.

\subsection{Entire Method}
For the reader's convenience, in this subsection we state our method in self-contained form in Algorithm~\ref{algo:Proposed}. Below, we provide an analysis for this algorithm, i.e. we show its time-complexity in terms of power iterations, and provide proof for the acclaimed probability at which the algorithm returns the correct result.

\begin{algorithm}
	\caption{Proposed method}
	\label{algo:Proposed}
	\begin{algorithmic}[1]
		\Procedure{Solver}{$\bA$}
			\State \textit{// - - - Bound of spectrum - - - }
			\State $k := \lceil 2.9 \cdot \log_{10}(N)+54.1\rceil$
			\State $\check{\mu}^2 := 0.5 / $\textsc{PowerIteration}($\bA^{-2},B_N,k$)
			\State $\hat{\mu}^2 := 1.5 \cdot $\textsc{PowerIteration}($\bA^2,B_N,k$)
			\State {It holds with probability $>1-10^{-9}$: } $\sigma(\bA)\subset(\check{\mu},\hat{\mu})$ and $\hat{\mu}/\check{\mu}\leq 3  \cdot \kappa(\bA)$\,.
			\State \textit{// - - - Bisection - - - }
			\State $J := \lceil\log_2(\hat{\mu}/\check{\mu})\rceil$
			\State $k := \lceil\,2.9 \cdot \log_{10}(N) + 1.8 \cdot \log_{2}(J) + 40.3\rceil$
			\For{$j=0,1,2,\dots,(J-1)$}
				\State $\mu := 1+1/2^j$
				\State $\bZ := (\mu-1) \cdot \big(\,(\mu-1) \cdot \bI + \frac{1}{\hat{\mu}}\cdot\bA\,\big)\inv$
				\State $[\tlambda,\bb] := $\textsc{PowerIteration}($\bZ,B_N,k$)
				\If{$\tlambda\geq 1$}
					\State \Return Matrix $\bA$ is not positive definite. Proof: $\bb\t\cdot\bA\cdot\bb\leq 0\,.$
				\EndIf
			\EndFor
			\State \Return Matrix $\bA$ is positive definite with probability $>1-10^{-6}$.
		\EndProcedure
	\end{algorithmic}
\end{algorithm}

\paragraph{Correctness and probability of the algorithm}
We discuss the probabilities in Algorithm~\ref{algo:Proposed} line by line. Therein, we always use $\varepsilon:=0.4$. Given that, we can simplify the formula for $k$ from \eqref{eqn:FormulaK} to the upper bound
\begin{align}
	k \leq \lceil 2.9 \cdot \log_{10}(N) + 5.76 \cdot \log_{10}(1/\delta)\rceil\,.\label{eqn:BoundK}
\end{align}

\largeparbreak

In lines 4 and 5 we compute an $\varepsilon$-accurate approximation of $\lambda_N(\bA),\blambda_1(\bA)$ with a probability $\delta$ of failure
\begin{align*}
	\delta_\kappa := 4 \cdot 10^{-10}\,.
\end{align*}
Hence, we obtain for $k$ from \eqref{eqn:BoundK} the suitable value
\begin{align*}
	k_\kappa := \lceil2.9 \cdot \log_{10}(N) + 54.1\rceil\,.
\end{align*}
According to Theorem~\ref{thm:PowerIter} and product rule of conditional probability, the probability of failure for the result in line 9 is
\begin{align*}
	\delta_0 := 1-(1-\delta_\kappa)^2 < 10^{-9}\,.
\end{align*}

\largeparbreak

In the next part we apply the bisection. According to the result in line 6, the value of $J$ satisfies
$$ 	J \in \cO\Big(\,\log\big(\kappa(\bA)\big)\,\Big) 	$$
with probability $>1-10^{-9}$.

The value of $k$ in line 9 is chosen from \eqref{eqn:BoundK} for a particular value $\delta = \tilde{\delta}$, that we derive now.

In the beginning of the $j$th for-loop in line 11, the bisection assumption
\begin{align}
	\lambda_{\text{max}}(\bM) \leq \mu_j \label{eqn:ProbJ}
\end{align}
holds with the conditional probability, that all former bisection iterations held true. Let us denote the probability, that \eqref{eqn:ProbJ} holds true, with $q_j$. Using the product rule of conditional probability, we find the following formula for $q_j$:
\begin{align*}
	q_j = (1-\delta_0) \cdot \tilde{q}^j
\end{align*}
In this formula, $\tilde{q}$ is the probability that the power iteration in line 13 provides an $\varepsilon$-accurate estimate for the eigenvalue. As we find from Theorem~\ref{thm:PowerIter}, this probability is
$$ 	\tilde{q} = 1-\tilde{\delta}\,. 	$$
The returned statement in line 18 holds true if and only if all power-iterations provided $\varepsilon$-accurate estimates for the largest eigenvalues. Hence, the probability that the statement in line 18 holds is
\begin{align*}
	q_{J} = (1-\delta_0) \cdot \tilde{q}^{J}\,.
\end{align*}
We want to propose an algorithm that is highly unlikely to provide a false result. Hence, for $q_J$ we require the following very high probability:
\begin{align*}
	q_J := 1-10^{-6}\,.
\end{align*}
Based on that, we obtain that $\tilde{\delta}$ is sufficiently small if it satisfies
\begin{align*}
	1-10^{-6} = (1-10^{-9}) \cdot (1-\tilde{\delta})^J\,.
\end{align*}
Defining 
$$ 	\rho := \frac{1-10^{-6}}{1-10^{-9}}\,,	$$
this can be reformulated as
\begin{align*}
	\tilde{\delta} = 1 - \rho^{(1/J)}\,.
\end{align*}
We now state the formula for $k$ in line 9:
\begin{align*}
	k = \lceil 2.9 \cdot \log_{10}(N) + 5.76 \cdot \log_{10}(1/\tilde{\delta})\rceil
\end{align*}
In the following we provide a bound for the right side of the equation. I.e., we provide an unnecessarily large value for $k$.

We can use $J\geq 1$ (because otherwise there would be no bisection iteration anyway) and hence obtain the bound
\begin{align*}
	\log_{10}(1/\tilde{\delta})
	&= -\log_{10}(\tilde{\delta})\\
	&= -\log_{10}\Big(\,1-\rho^{(1/J)}\,\Big) < 7 + \log_{10}(J)\,.
\end{align*}
All in all, a suitable value for $k$ is found as
\begin{align*}
	k 	&= \lceil\,2.9 \cdot \log_{10}(N) + 5.76 \cdot 7 + 5.76 \cdot \log_{10}(J)\\
		&< \lceil\,2.9 \cdot \log_{10}(N) + 1.8 \cdot \log_{2}(J) + 40.3\rceil =: \tilde{k}\,.
\end{align*}
In conclusion, the overall probability of correctness of the proposed method, with the particular formulas for $k$ in lines 3 and 9, is $>1-10^{-6}$.

\paragraph{Number of power iterations}
We now comment on the time-complexity. We count the complexity as the number $N_\text{for}$ of total power iterations with $\bA$ and the number $N_\text{inv}$ of total power iterations with a shifted inverse of $\bA$.

The number of matrix-vector products with $\bA$ equals $2 \cdot k_\kappa$, the value of $k$ in line $3$. It lives in the complexity
\begin{align*}
	N_\text{for} \in \cO\big(\,\log(N)\,\big)\,.
\end{align*}
In practice, this number is roughly bounded as
$$ 	N_\text{for} \lesssim 6 \cdot \log_{10}(N) + 108\,. 	$$

The number of power iterations with a shifted inverse of $\bA$ is $2 \cdot k_\kappa + J \cdot \tilde{k}$\,. This number has the complexity
$$ 	N_\text{inv} \in \cO\Big(\,\log(N) \cdot \log\big(\kappa(\bA)\big)^2 \cdot \log\left(\log\big(1+\kappa(\bA)\big)\right)\,\Big) 	$$
with a probability of $>1-10^{-9}$. The probability comes into play because there is an unlikely possibility that $J$ differs hugely in magnitude from $\kappa(\bA)$. If we assume $\kappa(\bA)\leq 10^{20}$, then in practice we obtain values like
$$ N_\text{inv} \lesssim \log_2\big(\kappa(\bA)\big) \cdot \big(\, 6 \cdot \log_{10}(N) + 108 \,\big)\,.	$$

\section{Conclusions}

We presented a new method for the numerical determination of positive definiteness of a symmetric matrix.
The algorithm either returns a proof, that the matrix is not positive definite, or it says that the matrix is positive definite. In the latter case, the probability, that the method's answer is wrong, is very small.

The time-complexity of the method is probabilistic in theory. But that is so mild and all damped by logarithms, that it can be ignored in practical terms. The cost of the method can be well-estimated as solving
$$ 	\log_2\big(\kappa(\bA)\big) \cdot \Big(\,6 \cdot \log_{10}(N) +100\,\Big) 	$$
linear equation systems with a shifted matrix of $\bA$.
While on first glance this does not appear very attractive, it can be quite efficient when $\bA$ is very large and when shifted linear systems with $\bA$ can be solved rapidly through parallel computations.

\paragraph{Promising Directions and Open Questions}
The new method forms a promising candidate for inertia correction of the linear systems that arise in optimization problems from discretized optimal control. In particular, we use the direct transcription approach in \cite{myOCP} to solve optimal control problems. There result non-linear programming problems, which are solved using Quasi-Newton based approaches. These form a large sparse saddle-point linear system, which in turn is reduced into a large thin-banded symmetric linear system. For the Newton step to yield descent, one must determine (and repair) positive definiteness of this large thin-banded linear system. The dimension $N$ of this matrix is in the order of $10^9$ to $10^{10}$. 

(Shifted) linear equation systems with this matrix can be solved rapidly, using a special parallel scalable version of the SPIKE algorithm \cite{myblocktridiagLGS}. In contrast to that, it would be entirely too time-consuming to perform a Choleskly decomposition, especially when it is used only to determine whether the system has the correct inertia. 

With the new method described in this paper, there is a chance that positive definiteness of the large thin-banded system can be determined on massively parallel machines in an affordable amount of time. Though, this of course will require a thorough revision of the parameters.

\largeparbreak

Open questions are related to finding a useful compromise for the complexity factors of our method and the probability $\delta$ of failure. In the above proposed form, we have chosen very small values for $\delta$. Certainly, cheaper variants of our method can be found by giving up on these restrictions and choosing a larger value $\delta$. Clearly, the best compromise between $\delta$ and the computational cost depends on the application at hand.

For our use case of inertia correction in Newton-type methods for constrained optimization, we usually do not bother too much if we add a shift to the Hessian that is actually larger than required. Hence, particularly for our applications we will perform further work towards a cheaper variant.

Further work in this direction may also be related to some initial heuristics, that can quickly solve those problems where $\bA$ has a spectrum that is advantageous for a cheap determination.

\FloatBarrier

\bibliography{PositiveDefiniteness_bib}
\bibliographystyle{plain}

\end{document}